\documentclass[12pt]{article}
cm \voffset -0.8 true cm

\usepackage{latexsym}
\usepackage{epic,eepic,epsfig}
\usepackage{amssymb}

\usepackage{amsfonts,srcltx}
\newtheorem{theorem}{Theorem}[section]
\newtheorem{prop}[theorem]{Proposition}
\newtheorem{lem}[theorem]{Lemma}
\newtheorem{rem}[theorem]{Remark}

\newtheorem{exam}[theorem]{Example}

\newtheorem{defi}[theorem]{Definition}

\newcommand{\proof}{{\sc Proof. \ }}

\newcommand{\mod}{\, \hbox{mod} \, }

\newcommand{\Aut}{\hbox{Aut}\, }

\newcommand{\Exp}{\hbox{Exp}\, }

\newcommand{\qed}{\hfill \mbox{\raisebox{0.7ex}{\fbox{}}}}

\def\demo{{\bf Proof}\hskip10pt}
\def\a{\alpha} \def\b{\beta} \def\g{\gamma} \def\d{\delta} 
 \def\s{\sigma} \def\t{\tau} \def\om{\omega} 
\def\th{\theta}   
  
\def\si{\Sigma} \def\O{\Omega}
\def\G{\Gamma}    
 \def\og{\overline G} \def\oh{\overline H}
 
  \def\om{\overline M}
\def\on{\overline N}

\def\o{\overline}

\def\o{\overline}

\def\of{\overline F} \def\od{\overline D}

\def\mb{\mathcal{B}}

\def\di{\bigm|} \def\lg{\langle} \def\rg{\rangle} 
\def\nd{\mathrel{\bigm|\kern-.7em/}}

\def\f{\noindent}

\def\PSL{\hbox{\em PSL}}
\def\PSU{\hbox{\em PSU}}

\def\P\GL{\hbox{\em P\GL}}
\def\SL{\hbox{\em SL}} \def\Exp{\hbox{\em Exp}}

\def\Aut{\hbox{\rm Aut}}

\def\Syl{\hbox{\rm Syl}}

\def\soc{\hbox{\rm soc}}

\def\Exp{\hbox{\rm Exp }}

\def\mod{\hbox{\rm mod }}
\def\Core{\hbox{\rm Core}}

\def\AG{\hbox{\rm AG}}
\def\PGL{\hbox{\rm PGL}}
\def\PSL{\hbox{\rm PSL}}\def\PSU{\hbox{\rm PSU}}
\def\AGL{\hbox{\rm AGL}}\def\SL{\hbox{\rm SL}}
\def\GL{\hbox{\rm GL}} \def\PG{\hbox{\rm PG}}
\def\GF{\hbox{\rm GF}}
\def\soc{\hbox {\rm soc}}

\begin{document}


\begin{center}
{\bf\Large SEMISYMMETRIC GRAPHS OF ORDER $2p^3$}
\end{center}

\begin{center}
{\sc Li Wang} and {\sc Shaofei Du}
\end{center}

\begin{center} {\small  School of Mathematical Sciences

Capital Normal University

Beijing, 100048, P R China}
\end{center}

\begin{abstract}
 A simple undirected graph is said to be {\em semisymmetric} if it is
regular and edge-transitive but not vertex-transitive. Every
semisymmetric graph is a bipartite graph with two parts of equal
size. It was   proved in  [{\em J. Combin. Theory
Ser. B} {\bf 3}(1967), 215-232] that there exist no semisymmetric graphs of order $2p$ and $2p^2$, where  $p$ is a prime.
 The classification of  semisymmetric graphs of order $2pq$  was given in [{\em Comm. in Algebra} {\bf 28}(2000), 2685-2715], for  any  distinct primes  $p$ and $q$.
  Our long term goal is to determine
  all the semisymmetric  graphs of order $2p^3$, for any prime $p$. All these graphs $\G$ are divided into two subclasses:  (I)  $\Aut(\G)$ acts
unfaithfully on at least one bipart; and  (II)  $\Aut(\G)$ acts
faithfully on both biparts.  This paper gives a group theoretical characterization  for  Subclass (I) and based  on this characterization, we   shall
give  a complete  classification
for this subclass  in our further research.
\end{abstract}


\section{Introduction}
All graphs considered in this paper are finite, undirected and
simple.  For a graph $\G$ with the  vertex set ${V}$ and  edge set $E$, by $\{u,v\}$ and $(u,v)$
we denote an edge and arc of $\G$, respectively, by $\Aut (\G)$ to denote its full automorphism group.  Set $A=\Aut(\G)$.
If $\G$ is bipartite with the bipartition   $V=W\cup U$, then we let $A^+$ be a subgroup of $A$ preserving  both $W$ and
$U$. Clearly if $\G$ is connected, then either $|A:A^+|=2$ or
$A=A^+$, depending on whether or not there exists an automorphism
which interchanges the two biparts.
For $G\le A^+,$ the graph $\G$ is said to be
$G$-{\em semitransitive} if $G$ acts transitively on both $W$ and
$U$, while  an  $A^+$-semitransitive graph is simply said  to be {\em semitransitive}.

A graph is said to be  {\em semisymmetric} if it is regular and
edge-transitive but not vertex-transitive.
It is easy to see that
every semisymmetric graph is a semitransitive bipartite graph with two biparts of
equal size.

The first person who studied semisymmetric graphs was Folkman. In
1967 he constructed several infinite families of such graphs and
proposed eight open problems (see \cite{Fol}). Afterwards, Bouwer,
Titov, Klin, I.V.~Ivanov, A.A.~Ivanov and others did much work on
semisymmetric graphs (see \cite{Bou1,Bou2,II,Iva,Kli,Tit}). They
gave new constructions of such graphs and nearly solved all of
Folkman's open problems. In particular, by using group-theoretical methods, Iofinova and Ivanov
\cite{II} in 1985 classified  cubic semisymmetric  graphs whose automorphism group acts primitively on both biparts, which
  was the first classification
theorem for semisymmetric graphs. More recently, following some deep results
in group theory which depend on the classification of finite simple
groups and some methods from graph coverings, some new  results of
semisymmetric graphs have appeared. For instance, in \cite{DX}  the second author
and Xu classified semisymmetric graphs of order $2pq$ for two
distinct  primes $p$ and $q$. For more results on semisymmetric
graphs, see
(\cite{CMMP,Du1,DM1,DM2,DWZ,DX,FK,LWX,LV,MMW,MMP1,Par,Wil} and so
on).

In \cite{Fol}, Folkman proved that there are no semisymmetric graphs
of order $2p$ and $2p^2$ where $p$ is a prime.  Then  we are
interested in determining semisymmetric graphs of order $2p^3$,
where $p$ is a prime. Since the smallest  semisymmetric graphs have order $20$ (see \cite{Fol}), we let $p\ge 3$.
  It was proved in \cite{MMW} that the Gray graph
of order $54$ is the only cubic semisymmetric graph of order
$2p^3$. To classify all the semisymmetric  graphs of order $2p^3$ is
still one of attractive and   difficult problems. These graphs $\G$
are naturally divided into two  subclasses:

 \begin{enumerate}
\item[{\rm }] Subclass (I):  $\Aut(\G)$ acts
unfaithfully on at least one bipart;
\item[{\rm }]  Subclass (II):  $\Aut(\G)$ acts faithfully on
both  biparts.
\end{enumerate}
The aim of this paper is to give   a group theoretical  characterization  for Subclass (I).  Based on  this characterization, we shall  give  a complete  classification
for this subclass  in our further research.

\vskip 3mm
In the following two paragraphs,  we first introduce two definitions used later.

Let $\cal P$ be  a partition  of the vertex set $V$.  Then we let $\G_{\cal P}$ be the {\it quotient graph} of $\G$ relative to $\cal P$, that is,
the graph with the vertex set $\cal P$, where two subsets $V_1$ and ${V_2}$ in $\cal P$ are adjacent  if there  exist  two vertices $v_1\in V_1$ and $v_2\in V_2$ such that $v_1$ and $v_2$ are adjacent in $\G$. In particular,  when $\cal P$ is  the set of orbits of a subgroup $N$ of $\Aut(\G)$, we denote $\G_{\cal P}$   by $\G_N$.

Let $\si =({\cal V}, {\cal E})$ be  a connected semitransitive and  edge-transitive graph  with
bipartition ${\cal V}={\cal W}\cup {\cal U}$, where $|{\cal W}|=p^3$ and
$|{\cal U}|=p^2$ for an odd   prime $p$.
Now we define a bipartite graph $\G =(V, E)$ with bipartition
 $V=W\cup U$, where
$$\begin{array}{ll}  &W={\cal W},\quad  U={\cal U}\times Z_p=\{ ({\bf u}, i)\di {\bf u}\in {\cal U}, i\in Z_p \},
\\&E=\{ \{ {\bf w}, ({\bf u}, i)\}\di \{ {\bf w}, {\bf u}\} \in {\cal  E}, i\in Z_p\}.$$
\end{array}$$
Then  we shall call that   $\G$   is  the {\it
graph  expanded from $\si $}.
Clearly $\G$ is
edge-transitive and regular.  By the definion, we see that  for any ${\bf u}\in {\cal U}$, the  $p$ vertices $\{ ({\bf u}, i)\di i\in Z_p\}$ in
$U$  have the same neighborhood in $\G $. Therefore,
$\G $  is semisymmetric,  provided there exist no two vertices
in ${\cal W}$ which have the  same neighborhood in $\si $.

\vskip 3mm
To state our main theorem, we first define four graphs: $\si (3)$, $\si (9)$, $\G(9)$ and $\G(18)$.
\begin{exam} \label{defi1} Let $\mathbb{V}=\mathbb{V}(3,3)$ be the 3-dimensional vector space over $\GF(3)$. Take three 2-dimensional subspaces of $\mathbb{V}:$
$$\begin{array}{lll}\mathbb{V}_0&=&\{ (0, b, c)\di b, c\in \GF(3)\}, \,  \mathbb{V}_1=\{ (a, 0, c)\di a, c\in \GF(3)\},\\
 \mathbb{V}_2&=&\{ (a, b, 0)\di a, b\in \GF(3)\} . \end{array}$$
Let  ${\cal W}=\mathbb{V}$ and let
${\cal U}=\{ \a +\mathbb{V}_i\di \a \in \mathbb{V}, i\in Z_3\}$, the set of  nine 2-dimensional subspaces (not all) in the 3-dimensional affine geometry $\AG(3,3)$ over $\GF(3)$.

 \vskip 3mm Define a  bipartite graph $\si (3)$  with biparts
${\cal W}$ and  ${\cal U}$, whose edge-set is
 $$\{ \{\a ,  \a +\mathbb{V}_i\}\di \a \in \mathbb{V}, i\in Z_3\} .$$

Define  $\si (6)$  to  be the bi-complement of $\si (3).$

Define  $\G (9)$ and $\G (18)$  to  be the  graphs expanded from $\si (3)$ and $\si(6)$, respectively.
\end{exam}

\begin{lem} \label{L} Both $\G(9)$ and $\G(18)$ are semisymmetric graphs of order 54, with   valency 9 and 18, respectively.
\end{lem}
\demo
Let $N$ be the translation group of the affine group $\AGL(3,3)$ and let   $L$ be the subgroup of $\GL(3,3)$ consisting of all  those $3\times 3$
  matrices with only one nonzero entry in each row and column.
Then it is easy to verify that   $N\rtimes L$ preserves the edge-set of both graphs   $\si (3)$ and $\si(6)$ and acts   edge-transitively on them.

Clearly, $\G(9)$ and $\G(18)$ are edge-transitive graphs  of 54, with  valency  9 and 18, respectively.
Since there exist no two vertices in ${\cal W}$ having the same neighborhood in $\si $ and since for each $i$, three vertices  $\{ (\a+\mathbb{V}_{i}, j)\di
j\in Z_3\}$ have the same neighborhood in $\G$, they are  not vertex-transitive and then  semisymmetric.\qed

\vskip 3mm

The main results  of this paper are the following Theorem~\ref{main1} and Theorem~\ref{main2}.

\begin{theorem} \label{main1}
 For  any odd prime $p$, let  $\G =(V, E)$ be  a semisymmetric graph of order $2p^3$  with the partition $V=W\cup U$ and full automorphism group  $A=\Aut(\G)$. Suppose that  $A$ acts unfaithfully on at least one bipart, say  $W$, with    the kernel  $A_{(W)}$. Let ${\cal W}=W$ and  $\cal U$  the  set of  orbits of $A_{(W)}$ on $U$. Set $\si=\G_{A_{(W)}},$   the quotient graph  of $\G$ induced by $A_{(W)}$, with the partition  ${\cal W}\cup {\cal U}$.
 Then the following hold.
 \begin{enumerate}
\item[{\rm (1)}] Every orbit of  $A_{(W)}$ on $U$ has length $p$,  $A_{(W)}\cong (S_p)^{p^2}$ and   $\G $ is expanded from $\si $.
 \item[{\rm (2)}]   $A/A_{(W)}$ acts faithfully on  ${\cal U}$ and so on ${\cal W}\cup {\cal U},$ and $A/A_{(W)}\cong \Aut(\si )$.
\item[{\rm (3)}]  $A$ acts faithfully  on $U$ and  there exist no two  vertices in $W$  having the same neighborhood in $\G $.
\end{enumerate}
\end{theorem}

By Theorem~\ref{main1}, we know that the graph $\G$  is uniquely determined by its quotient graph $\si $.
 Now  we turn to focus on the graph $\si $, see    the following  theorem.
\begin{theorem} \label{main2} Adopting  the notation in Theorem~\ref{main1} and setting $F=\Aut(\si)$, we have
 \begin{enumerate}
\item[{\rm (1)}] $F$ acts   imprimitively on $\cal U$.
\item[{\rm (2)}] Suppose that $F$  acts   primitively on $\cal W$.  Then $p=3$, and $\si \cong \si(3)$ or $\si(6)$;  $\G \cong \G(9)$ or  $\G(18)$, see Example~\ref{defi1}.
\item[{\rm (3)}] Suppose that    $F$  acts imprimitively on  $\cal W$, with a block (of length $p$) system $\mathfrak{U}$ and the kernel $F_{(\mathfrak{U})}$.
 Then   either
  \begin{enumerate}
\item[{\rm (3.1)}] $F_{(\mathfrak{U})}$  is solvable and  acts transitively on $\cal W$ and  $F$ is an affine group;  or

\item[{\rm (3.2)}] $F_{(\mathfrak{U})}$  induces blocks of length $p^2$ on  $\cal W$.
 Take  any  $\bf w$ in $\cal W$. Then
 either
\begin{enumerate}
\item[{\rm (3.2.1)}] $\bf w$ is exactly adjacent to two blocks in $\mathfrak{U}$; or
\item[{\rm (3.2.2)}]  $\bf w$ is adjacent to at least three blocks in $\mathfrak{U}$,  $p\ge 5$,  $F$   contains the nonabelian normal  $p$-subgroup acting regularly  on $\cal W$, $F_{(\mathfrak{U})}$ is solvable, and $F/F_{(\mathfrak{U})}\cong Z_p\rtimes Z_r$, where $r\di (p-1)$ and $r\ge 3$.
\end{enumerate}
\end{enumerate}
\end{enumerate}
\end{theorem}

\begin{rem} To classify  all the graphs $\G$  in  Subclass (I), it suffices to
 determine the   graphs $\si $ in  Theorem~\ref{main2}.(3.1), (3.2.1) and (3.2.2). However, the determination of these graphs  is still quite complicated. In this paper, we  just construct  the respective examples, see    Section 5, and  by using the group structures obtained  in Theorem~\ref{main2},  we shall   give  a   complete classification for  them  in our further research.
\end{rem}

After this introductory section,
some preliminary results will be given in Section 2;  Theorem~\ref{main1} and Theorem~\ref{main2} will be proved in Section 3 and
 4, respectively. Finally, the related  graphs will be constructed in Section 5.

\section{Preliminaries}
 First we introduce some notation:
 by $K_n$ and $K_{m,n}$  we denote the complete graph of order $n$ and the complete bipartite graph with two biparts of size $m$ and $n$, respectively.
 For a graph $\G$,  by $d(v)$ we denote the degree of a vertex $v\in V$.
 For a prime $p$, by $p^i\di\di n$ we mean
$p^i\di n$ but $p^{i+1}\nmid n$. By
$Z_n$, $D_{2n}$ and $S_n$, we denote the cyclic group of order $n$, the dihedral group of
order $2n$  and the symmetric group of degree $n$, respectively. By $\GF(p)$, we denote the field of $p$ elements. For a ring $S,$  let $S^*$ be the
multiplicative group of all the units in $S.$
For  a transitive group $G$  on $\O $ and  a subset $\O_1$ of $\O$,  by
 $G_{\O_1}$ and $G_{(\O_1)}$ we denote the setwise  stabilizer and pointwise   stabilizer of $G$ relative to $\O_1$, respectively.
 A $m$-block of $G$ means a block with length $m$.

For a group  $G$ and  a subgroup $H$ of $G$,
use $Z(G),$ $C_G(H)$ and $N_G(H)$ to denote the center of $G$, the
centralizer and normalizer of $H$ in  $G,$ respectively.
  A semidirect product of the group $N$ by the group $H$ is denoted by $N\rtimes H,$ where $N$ is normal.
  A wreath product of $N$ by $H$ is denoted by $N\wr H$, that is $N^n\rtimes H$, where $H\le S_n$.
By $[G:H]$ we
denote the set of right cosets of $H$ in $G$. The action of $G$ on
$[G:H]$ is always assumed to be the right multiplication action.

For any $\a$ in the $n$-dimensional  vector space  $\mathbb{V=}\mathbb{V}(n,p)$ over  $\GF(P)$, we denote by $t_\a $
the translation corresponding to $\a $ in the affine geometry $\AG(\mathbb{V})$  and  by $T$ the
translation subgroup of the affine group $\AGL(n,p)$. Then $\AGL(n, p)\cong T\rtimes
\GL(n, p).$ We adopt matrix notation for $\GL(n, p)$ and so
 we have
$g^{-1}t_\a g=(t_\a)^g=t_{\a g}$ for any $t_\a \in T\le \AGL(n,p)$
 and  $g\in \GL (n, p).$

For group-theoretic concepts and notation not defined here the
reader is refereed to \cite{Dix,Hup}.

\vskip 3mm To  constructed    graphs,  we need to introduce  the definition of bi-coset graphs and
 two  properties.
\begin{defi} \label{defi2} \cite{DX}
Let $G$ be a group,   $L$ and $R$    subgroups of $G$ and let
$D=RdL$ be a  double coset of $R$ and $L$ in $G$. Let
$[G:L]$ and $[G:R]$ denote the set of right cosets of $G$ relative to $L$ and $R$ respectively.
  Define a bipartite graph $\G={\bf B}(G,L,R;D)$ with
bipartition $V=[G:L]\cup [G:R]$ and  edge set
$E=\{ \{ Lg,Rdg\}\di g\in G, d\in D\}$.
 This graph is called the
bi-coset graph of $G$ with respect to $L$, $R$ and $D$.
\end{defi}

\begin{prop} \label{property}\cite{DX} 
The graph $\G={\bf B}(G,L,R;D)$ is a well-defined bipartite graph.
Under the right multiplication action on $V$ of $G$, the graph
$\G$ is $G$-semitransitive.  The kernel of the action of $G$ on
$V$ is $\Core_G(L)\cap\Core_G(R)$, the intersection of the cores
of the subgroups $L$ and $R$ in $G$. Furthermore, we have

\begin{enumerate}
\item[{\rm (i)}] $\G$ is $G$-edge-transitive;
\item[{\rm (ii)}]  the degree of any vertex in $[G:L]$ (resp. $[G:R])$
is equal to the number of right cosets of $R$ (resp. $L$) in $D$
(resp. $D^{-1})$, so $\G$ is regular if and only if $|L|=|R|;$
\item[{\rm (iii)}]  $\G$ is connected if and only if $G$ is generated
by elements of $D^{-1}D$.
\end{enumerate}
\end{prop}

\begin{prop}\label{iden} \cite{DX}
Suppose $\G'$ is a $G$-semitransitive and edge-transitive graph with bipartition
$V=U$ $\cup\, W$. Take $u\in U$ and $w\in W$. Set \
$D=\{ g\in G\di w^g \in \G'_1(u)\}.$ \ Then $D=G_wgG_u$ and $\G'\cong {\bf B}(G, G_u, G_w;D).$
\end{prop}

\vskip 3mm
Finally, several group theoretical results are given.

\begin{prop}\cite{Gur}\label{primep}
   Let T be a nonabelian simple group with a subgroup $H<T$ satisfying
 $|T:H|=p^a,$ for $p$ a prime. Then one of the following holds:
\begin{enumerate}\item[{\rm (i)}]  $T=A_n$ and $H=A_{n-1}$ with
$n=p^a;$
\item[{\rm (ii)}]  $T=\PSL(n,q)$, $H$ is the stabilizer
of a projective point or a hyperplane in $\PG(n-1,q)$ and
$|T:H|=(q^n-1)/(q-1)=p^a;$ \item[{\rm (iii)}]  $T=\PSL(2,11)$ and
$H=A_5;$ \item[{\rm (iv)}]  $T=M_{11}$ and $H=M_{10};$ \item[{\rm
(v)}]  $T=M_{23}$ and $H=M_{22};$ \item[{\rm (vi)}] $T=\PSU(4,2)$
 and $H$ is a subgroup of index 27.\end{enumerate}
\end{prop}

\begin{prop}\cite{Blo}\label{PSL3p}
For an odd prime $p$, let $\og =\PSL(3,p)$ and  $\oh $  a  proper subgroup of $\og
$. Then one of the following holds:

\begin{enumerate}
\item[{\rm (I)}]  If $\oh $ has no nontrivial normal elementary
abelian subgroup, then $\oh $ is conjugate in $\GL(3,p)/Z(\SL
(3,p))$ to one of the following groups:

\begin{enumerate}
\item[{\rm (i)}]  $\PSL(2,7)$, with $p^3\equiv 1 (\mod  7);$
\item[{\rm (ii)}]  $A_6$, with $p\equiv 1, 19(\mod  30);$
\item[{\rm (iii)}] $\PSL(2,5)$, with $p\equiv \pm 1 (\mod  10);$
 \item[{\rm (iv)}]  $\PSL(2, p)$ or $\PGL(2,p)$  for $p\ge 5$.
\end{enumerate}

\item[{\rm (II)}]   If $\oh $ has a nontrivial normal elementary
abelian subgroup, then
 $\oh $ is conjugate to a subgroup of one of the following subgroups:

\begin{enumerate}
\item[{\rm (i)}]  $Z_{(p^2+p+1)/(3, p-1)}\rtimes Z_3;$
 \item[{\rm (ii)}] the subgroup $\of $ of all matrices with only one nonzero
entry in each row and column, and $\of $ contains the subgroup
$\od $ of all diagonal matrices as a normal subgroup such that
$\of /\od \cong S_3;$
 \item[{\rm (iii)}] the point- or line-stabilizer of a given point
  $\lg (1,0,0)^T\rg $ or the line
$\lg (0, $ $\a, \b )^T\di \a, \b \in F_p\rg $;
 \item[{\rm (iv)}]  the group $\om $ such that $\om$ contains a normal subgroup $\on \cong Z_3^2$ and
 $\om /\on $ is isomorphic to $\SL(2,3)$ if $p\equiv 1(\mod 9)$ or to
 $Q_8$ if $p\equiv 4, 7(\mod 9)$.
 \end{enumerate}
\end{enumerate}\end{prop}

\begin{prop} \cite{DK}
\label{GL2p1} For an odd prime $p$, let $H$ be a maximal subgroup of
$G=\GL(2,p)$ and $H\ne \SL(2,p)$. Then up to conjugacy, $H$ is
isomorphic to one of the following subgroups:

\begin{enumerate}
 \item[{\rm (i)}] $D\rtimes \lg b\rg ;$
where $D$ is the subgroup of diagonal matrices and $b={\footnotesize
\left(\begin{array}{ll} 0 &1\\ 1 &0\end{array}\right)};$
 \item[{\rm (ii)}] $\lg a\rg \rtimes \lg b\rg $, where $b={\footnotesize
\left(\begin{array}{ll} 1 &0\\ 0 &-1\end{array}\right)}$ and
$\lg a\rg $ is the Singer subgroup of $G$, defined by
$a$={\footnotesize $\left(\begin{array}{ll} \g &\d\th \\ \d
&\g\end{array}\right)$}$\in G,$ where $F_p^*=\lg \th \rg $,
$F_{p^2}=F_p({\bf t})$ for ${\bf t} ^2=\th ,$ and  $F_{p^2}^*=\lg
\g+\d{\bf t}\rg;$
 \item[{\rm (iii)}]  $\lg a\rg
\rtimes D$, where $a={\footnotesize \left(
\begin{array}{cc} 1 &1 \\ 0&1 \end{array} \right)};$
 \item[{\rm (iv)}] $H/\lg z\rg $ is isomorphic to
 $A_4\times Z_{\frac{p-1}2}$,  for $p\equiv 5(\mod 8);$
  $S_4\times Z_{\frac{p-1}2}$ for  $p\equiv 1, 3, 7 (\mod 8);$
  or $A_5\times Z_{\frac{p-1}2}$ for $p\equiv \pm 1(\mod 10),$ where $z={\footnotesize
\left(\begin{array}{ll} -1 &0\\ 0 &-1\end{array}\right)}$,
  $Z_{\frac{p-1}2}=Z(G)/\lg z\rg ;$
\item[{\rm (v)}]  $H/\lg z\rg =A_4\rtimes \lg s\rg $, $\lg s^2\rg
\le Z(G)/\lg z\rg $, if $p\equiv 1(\mod 4).$
\end{enumerate}\end{prop}
The following theorem can be extracted from \cite{Dob}.
\begin{prop} \label{Dob}
\label{dobson} Let $G$ be a transitive permutation group of degree $p^2$, where $p\ge 5$ a prime and let $P$ be a Sylow $p$-subgroup of $G$.
Suppose that $G$ is imprimitive and $|P|=p^3$. Then $P\lhd G$.
\end{prop}

\section{Proof of Theorem~\ref{main1}}
From now on,   we assume that
 $\G$ is a semisymmetric graph of order $2p^3$ with the bipartition $V=W\cup U$, where $p$ is a prime, and  $A=\Aut(\G)$  acts unfaithfully on at least one part, say $W$.

 For avoiding confusions,    we  need to emphasis the following notation:

 $u:$  a  vertex in  $U$;    \,  ${\bf u}$:  a block induced by $A_{(W)}$ on $U$;

 ${\cal U}$: the  set of such blocks ${\bf u}$;

  \, $\mathfrak{u}$:  a block contained in ${\cal U}$; \,  $\mathfrak{U}$:  the  set of all such blocks $\mathfrak{u}$.

\vskip 3mm \f Symmetrically,  for other bipart $W$, we let $w,$  ${\bf w},$ ${\cal W},$  $\mathfrak{w}$ and  $\mathfrak{W}$ have the same meaning.
 Moreover, when emphasizing  on the set $\cal U$, we prefer to call ${\bf u}$ a {\it vertex  } in  $\cal U$ but not a {\it block} in $U$.
 \vskip 3mm

\vskip 3mm

\f {\bf Proof of Theorem~\ref{main1}:} Now  $A_{(W)}$ induces a complete $m$-block  system $${\cal U}=\{ {\bf u_0},  {\bf u_1}, \cdots,
 {\bf u_{\frac{p^3}m-1}}\},$$
  on $U,$ where  $m\di p^3$.  Let $A_{(\cal U)}$ be the kernel of $A$ on ${\cal U}$.
 Then  we  divide  the proof into  the following six steps.

  \vskip 3mm
 {\it Step 1:} Show that $p\ge 3$,  $\G\not\cong K_{p^3,p^3}$ and $\G$  is  connected.

 \vskip 3mm

 By ~\cite{Fol},  there exists no semisymmetric graph of  order less than 20.
Hence $p\ge 3$. Since the complete bipartite  graph $K_{p^3,p^3}$ is a symmetric graph,
$\G \not\cong K_{p^3,p^3}$.

Suppose that  $\G $  is disconnected. From the edge-transitivity of $\G $,  we get that $\G $ is isomorphic to either
$p^3K_2$ or $\frac {p^3}m\Phi_{2m}$ where $\Phi_{2m}$ (for $m\in \{ p,
p^2\} )$ is a regular edge-transitive bipartite graph of order $2m$. However,
 by  ~\cite{Fol}, there exist no semisymmetric graphs of order $p$ and $2p$, that is, $\Phi_{2m}$ is vertex-transitive,
 which implies  $\G$ is
vertex-transitive, a contradiction. Therefore,  $\G $  is  connected.

 \vskip 3mm
 {\it Step 2:} Show that $m\ne 1, p^3$.

 \vskip 3mm
 Since $A$ acts unfaithfully on $W$, we get $m\ne 1.$
 Suppose that  $m=p^3$.
Take $w\in W$. Since $A_{(W)}$ fixes $w$ and acts transitively
on $U$, it follows that $w$ is adjacent to all the vertices in $U$, which implies  $\G \cong K_{p^3,p^3},$ a contradiction.

\vskip 3mm
 {\it Step 3:} Show that $A_{(\cal U)}$ acts intransitively on $W$.
 \vskip 3mm
 Suppose that $A_{(\cal U)}$ acts transitively
on $W$. Then we shall show that for any $w\in W$, we have $\G_1(w)=U$, which implies  $\G\cong K_{p^3,p^3}$, a contradiction.

For any block ${\bf u}_j$ in ${\cal U}$, take an edge $\{w_1,u_j\}$ where $w_1\in W$ and $u_j\in {\bf u}_j$.
  Since $A_{(\cal U)}$ fixes ${\bf u}_j$ setwise and acts transitively
on $W$, there exists a $g\in A_{(\cal U)}$ sending $\{w_1,u_j\}$ to $\{w,u_j^g\}$, which means that $w$  is
adjacent to a vertex $u_j^g$ in ${\bf u}_j$.
Moreover,   since  $\{w,u_j^g\}\in E$ and  since $A_{(W)}$ fixes $w$ and acts transitively
on ${\bf u}_j$, it follows that $w$ is adjacent to all the  vertices in ${\bf u}_j$. Therefore, $\G_1(w)=U$.

\vskip 3mm
 {\it Step 4:} Show Theorem~\ref{main1}.(1).
\vskip 3mm
For the contrary,  suppose  $m=p^2.$  Then  $|{\cal U}|=p$  and
$A/A_{(\cal U)}\lessapprox S_p$. Moreover, since $A/A_{(W)}$ acts transitively on
$W$, we get $p^3\di |A/A_{(W)}|.$ Therefore, $A_{(W)}\lneqq A_{(\cal U)}$. Now
$A_{(\cal U)}$ induces a complete $n$-block system
${\cal W}=\{{\bf w}_1,\cdots, {\bf w}_{p^3/n}\}$ on $W.$ Since
$A/A_{(\cal U)}$ is transitive on $\cal W$ and  $A/A_{(\cal U)}\le
S_p$ , we get that $|{\cal W}|=p$.
 Then $n=p^2.$

 Consider the quotient  graph $\G_{A_{(\cal U)}}$ induced by  $A_{(\cal U)}$ with the  bipartition  ${\cal U} \bigcup {\cal W}.$
   Take the edge $\{{\bf u}_i, {\bf w}_j\}$ in $\G_{A_{(\cal U)}}$.
 A same argument as in  the proof of Step $1$
shows that the  induced subgraph $\G ({\bf u}_i\bigcup {\bf w}_j)\cong K_{p^2, p^2}.$ Therefore, for any ${\bf u}_{i_1}\in {\cal U}$ and ${\bf w}_{j_1}\in {\cal W}$,
  the induced subgraph $\G({\bf u}_{i_1}\bigcup {\bf w}_{j_1})$ is either an  empty graph or a complete bipartite
graph,  which implies   $\Gamma\cong \G_{A_{(\cal U)}}[p^2K_1].$  Since  the  graph $\G_{A_{(\cal U)}}$ of order $2p$ is  symmetric
 by~\cite{Fol}, we get $\Gamma$ is vertex-transitive, a contradiction again.

  Since $A_{(W)}$ fixes $W$ pointwise and acts transitively  on each ${\bf u}_i$ in ${\cal U}$,   it follows that  $p$ vertices  in each ${\bf u}_i$  have the same neighborhood in $\G$. Therefore, $\G $ is expanded from $\si $. Moreover,
  $A_{(W)}\cong S_p^{p^2}.$

\vskip 3mm
 {\it Step 5:} Show  Theorem~\ref{main1}.(2).

\vskip 3mm

From Step 4, we get $m=p$ and then $|{\cal U}|=p^2.$  Assume the contrary, that  is,
$A/A_{(W)}$ acts unfaithfully on $\cal U$.  Then $A_{(W)}\lneqq A_{(\cal U)}$.
As before, let ${\cal W}=\{{\bf w}_1,\cdots, {\bf w}_{p^3/n}\}$ be a complete
$n-$block system of $A_{(\cal U)}$ on $W.$

If $n=p,$ then $|{\cal W}|=p^2$ and as in Step 4 again, one may
easily see $\Gamma\cong \G_{A_{(\cal U)}}[pK_1],$ a contradiction.

Suppose  that  $n=p^2.$ Then $|{\cal W}|=p.$
 Then $A/A_{(\cal W)}\lessapprox
S_p$ and $A_{(\cal U)}\le A_{(\cal W)}$.
 Moreover, since $A/A_{(\cal U)}$ acts
transitively on ${\cal U}$, it follows that $p^2\di |A/A_{(\cal U)}|.$ Then
$A_{(\cal U)}\lneqq A_{(\cal W)}$. Naturally, we consider two cases:

\vskip 3mm (i) $A_{(\cal W)}$ is transitive on ${\cal U}$.

On the one hand, since $A_{(\cal W)}$ fixes $\cal W$ pointwise and acts
transitively on $\cal U$ and since $|{\cal W}|=p\ne |{\cal U}|=p^2$, the quotient graph of $\G$ with partition
${\cal W}\cup {\cal U}$ is isomorphic to $K_{p,p^2}$. On the other hand,
for any block ${\bf u}_i\in \cal U$ and ${\bf w}_j\in \cal W$, by
considering the actions of $A_{(W)}$ and $A_{(\cal U)}$ we know that the
induced subgraph $\G ({\bf u}_i\cup {\bf w}_j)$ is complete bipartite. Therefore,  $\Gamma \cong K_{p^3,p^3}$  a contradiction.

\vskip 3mm
 (ii) $A_{(\cal W)}$ has the blocks of length $p$ on ${\cal U}$.

Suppose that $A_{(\cal W)}$ has  blocks of length $p$ on ${\cal U}$. Then $A_{(\cal W)}$ has   blocks of length $p^2$ on $U$.
Then the quotient graph $\G_{A_{(\cal W)}}$ induced by $A_{(\cal W)}$ is  an edge-transitive graph of order $2p$ and then it
is symmetric  by~\cite{Fol} again.

Similarly, by considering the actions of $A_{(W)},$ $A_{(\cal U)}$ and
$A_{(\cal W)}$, we may show that the induced subgraph $\G ({\bf u}_i\cup
{\bf w}_j)$ is either complete bipartite or empty. Therefore, the graph
$\Gamma$ is vertex-transitive, a
contradiction.

 This proves that  $A/A_{(W)}$ acts faithfully on ${\cal U}.$

\vskip 3mm

Finally we  show that  $\Aut(\si)\cong A/A_{(W)}$.  Since  $A/A_{(W)}$  acts faithfully  on $\cal U$,
 it   induces a  faithful and  edge-transitive action on  $\si $, that is $A/A_{(W)}\lesssim \Aut(\si ).$
Clearly,   the graph $\G$ is uniquely determined by its the graph $\si$.
Then one  may see that every automorphism of $\si $ can be  extended to an automorphism of $\G $ which preserves ${\cal W}$, that means
$|\Aut(\si)|\le |A/A_{(W)}|$.
Therefore, $\Aut(\si)\cong A/A_{(W)}$.

\vskip 3mm
 {\it Step 6:} Show  Theorem~\ref{main1}.(3).

\vskip 3mm
 Since $A/A_{(W)}$
acts faithfully on ${\cal U}$ by Step 5, it follows that $A_{(\cal U)}=A_{(W)}$ and so
$A_{(U)}\le A_{(W)}$.  Since $A_{(U)}\cap A_{(W)}=1,$ we get $A_{(U)}=1$, equivalently, $A$
acts faithfully on $U$.

Suppose that there exist two vertices ${\bf w_1}$ and ${\bf w_2}$ in ${\cal W}$ having the same neighborhood in $\si $. Then the
 permutation $\t $  exchanging ${\bf w_1}$ and ${\bf w_2}$ and fixing
 other vertices of $\si $  is clearly an automorphism  of  $\si $, which
forces that $A/A_{(W)}$ acts unfaithfully on ${\cal U}$, a contradiction. Therefore, there exist no  two vertices ${\bf w_1}$ and ${\bf w_2}$ in ${\cal W}$ having the same neighborhood in $\si $.\qed

\section{Proof of Theorem~\ref{main2}}
By Theorem~\ref{main1},  from now on we focus on the quotient graph $\si $ induced by
 $A_{(W)}$ with biparts ${\cal W}\cup {\cal U}$, where $|{\cal W}|=p^3$ and $|{\cal U}|=p^2$.
 Since ${\bf w}=\{w\}$ for some $w\in W$, we shall   identify  ${\bf w}$ with $w$, and ${\cal W}$  with  $W$ as well.
Moreover,
 $\si $ is edge-transitive and there exist no  two vertices in ${\cal W}$ having the same neighborhood in $\si $.

To prove Theorem~\ref{main2}, we shall prove that
   $F=\Aut(\si)$ acts imprimitively on
${\cal U}$ in Subsection 4.1, that is Theorem~\ref{main2}.(1); and deal  with the cases when
   $\Aut(\si)$ acts primitively   on ${\cal W}$ in Subsection 4.2, that is Theorem~\ref{main2}.(2),  and imprimitively   on ${\cal W}$   in Subsection 4.3, that is Theorem~\ref{main2}.(3), respectively.

\subsection{Proof of Theorem~\ref{main2}.(1)}
First we prove a  group theoretical result.

\begin{lem}\label{p2-p3} For an odd prime  $p$, let $G$ be a primitive group on $\O$, where $|\O|=p^2$.  Suppose that $G$ has
a faithful transitive representation of degree $p^3$. Then $G$ is isomorphic to  one of the following  groups:
 \begin{enumerate}
\item[{\rm (1)}] ${\rm P}{\rm \G}{\rm L}(2,8)$, for $p=3$;
\item[{\rm (2)}] $Z_3^2\rtimes H$, where $H=\SL(2,3)$ or $\GL(2,3)$, for $p=3$;
\item[{\rm (3)}]  $Z_5^2\rtimes H$, where $H=\SL(2,5)$ or $\GL(2,5)$, for $p=5$;
\item[{\rm (4)}]  $Z_7^2\rtimes \SL(2,7)$, for $p=7$;
\item[{\rm (5)}]  $Z_{11}^2\rtimes \SL(2,11)$, for $p=11$.
\end{enumerate}
All these representations are imprimitive.\end{lem}
\demo
By the well-known O'Nan-Scott Theorem \cite{Dix}, every  primitive group $G$ of degree $p^2$  is
almost simple type, product type or affine type. Let $T=\soc(G).$
Suppose $G$ has a faithful transitive representation on $\O'$, where $|\O'|=p^3$.
Then we divided the proof into the following three cases.

\vskip 3mm
Case 1:   $G$ is almost simple type.

\vskip 3mm
In this case, $T=\soc(G)$ is either $A_{p^2}$ or $\PSL(n,q)$,
where $\frac{q^n-1}{q-1}=p^2$, by checking Proposition~\ref{primep}.
First suppose that $G$ is primitive on $\O'$. Then  by
checking  Proposition~\ref{primep} again, the  almost simple groups
of degree $p^3$ are: $A_{p^3}$, $\PSU(4,2)$ or
$\PSL(n,q),$ where $\frac{q^n-1}{q-1}=p^3$. Clearly, our group $G$ now
cannot have any faithful primitive representation of degree $p^3$.

In what follows, suppose that $G$ acts imprimitively on $\O'$. Let ${\cal B}$ be an imprimitive complete $m-$block system.
Then $T$ acts transitively on ${\cal B}$ with the kernel $K.$
 Since $T$ is the unique minimal normal subgroup of $G$, it follows that either $T\le K$
or $K=1.$ In other words, if $T$ acts transitively on $\O'$, then $K=1;$
if $T$ is intransitive on $\O'$, then $T\le K$ and  $p||G:T|$.

 (i)  Firstly, suppose that $T=A_{p^2}.$ Since $|G:T|\le 2,$  we can get that $T$ is impritimitive and transitive on $\O'.$
In this case, $K=1.$ Then $|{\cal B}|=p^2$ and $m=p.$
Take a block ${\bf b}$ in ${\cal B}.$ Then $T_{\bf b}=A_{p^2-1},$  which should be transitive on $\bf b$. However,  $A_{p^2-1}$ has no subgroup of index $p$, a contradiction.

(ii) Secondly, suppose that $T=\PSL(n,q)$, where $\frac{q^n-1}{q-1}=p^2$. Then $T\le G\le  {\rm P}\Gamma {\rm L}(n,q)$, where
$q=p_1^k$ for a prime $p_1.$  It is known that
$| {\rm P}\Gamma {\rm L}(n,q):T|\di (q-1)k$. If $n=2$, then $p^2-1=q=p_1^k$. From this
equation, we can get $p=3$ and $q=8,$ that is $T=\PSL(2,8)$ and
${\rm P}{\rm \G}{\rm L}(2,8)=T\rtimes \lg f\rg $, where $f$ is the field automorphism of order 3
of $T$. This is a case in  (1) of the lemma.

 Suppose that $n\ge 3$.  Then $p^2=\frac{q^n-1}{q-1}\ge q^2+q+1$ and so we
get $p>q=p_1^k$, which implies $p\nmid (q-1)$ and $p\nmid k$,
and then $p\nmid |G:T|$, that is,  $T\nleq K$ and thus $T$ acts transitively on $\O'.$
However, in what follows we shall show that   $p^3\nmid |T|$.

Since
$$|\PSL(n,q)|=\frac{(q^n-1)(q^n-q)\cdots(q^n-q^{n-1})}{(q-1)(n,q-1)}=p^2\frac{(q^n-q)\cdots(q^n-q^{n-1})}{(n,q-1)}$$ and
since $p\nmid q$ and $p\nmid (q-1),$
it suffices to show  $p\nmid (q^l+q^{l-1}+\cdots+1)$ for any $1\le l<n-1.$

Suppose that $k$ is the  minimal positive integer  such that $p\mid (q^k+q^{k-1}+\cdots+1).$
Write $$\begin{array}{lll}&&p^2=q^{n-1}+q^{n-2}+\cdots+1\\
&=&(1+q+q^2+\cdots+q^k)+q^{k+1}(1+q+q^2+\cdots+q^k)\\
&&+\cdots +q^{n-i}(1+q+q^2+\cdots+q^i).\end{array}$$
Then it follows $p\di (1+q+q^2+\cdots+q^i)$. From the minimality of $k$, we get   $i=k$ so that
$$p^2=(1+q+q^2+\cdots+q^k)(1+q^{k+1}+q^{2(k+1)}+\cdots+q^{n-k}),$$
and so
$$1+q+q^2+\cdots+q^k=1+q^{k+1}+q^{2(k+1)}+\cdots+q^{n-k},$$
that is
$$q(1+q+\cdots+q^{k-1}-q^{k}-q^{(2k-1)}-\cdots-q^{n-k-1})=1,$$
a contradiction.

\vskip 3mm
Case 2:   $G$ is product type.
\vskip 3mm
In this case, $G=(M\times M)\rtimes Z_2,$ where $M$ is an irregular primitive group of degree $p.$
Clearly,  $p^3\nmid |G|,$ a contradiction.

\vskip 3mm
Case 3:   $G$ is affine type.
\vskip 3mm
Now $G=N\rtimes H,$ where $N\cong Z_p^2$ and $H$ is an
irreducible subgroup of $\GL(2,p).$ Clearly, $G$ acts imprimitively on $\O'.$ Since $|\GL(2,p)|=p(p-1)^2(p+1),$
  we know that $p^3\di \di |G|$ and then $N$ induces  a $p^2-$block system, say  ${\cal B}$,  on $\O'$.
  Take ${\bf b}\in {\cal B}$ and $b\in {\bf b}$.   Considering
  the action of $H$ on ${\cal B}$, we know that $|H:H_{\bf b}|=p$ and then
  $H_{\bf b}=H_b$, that is, $H$ has a subgroup of index $p$. Checking Proposition~\ref{GL2p1}, we get that $H=\SL(2,p)$ for
   $p=3,5$, 7 and 11; or $H=\GL(2,p)$ for $p=3, 5$. This completes the proof of the lemma.

\qed

\vskip 3mm
\f {\bf Proof of Theorem~\ref{main2}.(1):}
For the contrary, suppose that $F$ acts primitively on ${\cal U}$.  Then $F$ has a
faithful primitive representation of degree $p^2$.  Since $|{\cal W}|=p^3$ ,  $F$ has a faithful  transitive
representation of degree $p^3$ and so $p^3\di |F|$. Then $F$ is one of the groups in   Lemma~\ref{p2-p3} and
we divide the proof into two cases according to   $F={\rm P}{\rm \G}{\rm L}(2,8)$ or $F$ is an affine group.

\vskip 3mm  (i) $F={\rm P}{\rm \G}{\rm L}(2,8)$

 \vskip 3mm
 Let  $F=T\rtimes \lg f\rg $ and let $H\cong Z_2^3\rtimes Z_7$ be a
point stabilizer of $T$ on the projective line. Then  $F_{\bf w}=H$ for some
${\bf w}\in {\cal W}$ and $F_{\bf u}=H\rtimes \lg f\rg $ for some ${\bf u}\in {\cal U}$. Since
each of $F_{\bf w}$, $F_{\bf w}^f$ and $F_{\bf w}^{f^2}$  fixes ${\bf u}$ and is  transitive on other 8 vertices on ${\cal U}$,
three vertices ${\bf w}$, ${\bf w}^f$ and ${\bf w}^{f^2}$ have the same neighborhood
in the graph $\si $, a contradiction (see Theorem~\ref{main1}.(3)).

\vskip 3mm
 (ii) $F$ is an affine group.

\vskip 3mm Now $F=N\rtimes H,$ where $N\cong Z_p^2$ and either $H=\SL(2,p)$ where $p=3, 5, 7, 11$ or $\GL(2,p)$ where $p=3,5 $.
First let  $H=\SL(2,p)$. Let $Z$ be  the center of
$\SL(2,p)$. Clearly, $Z\le H_{\bf w}$. If $p=3$ then $\frac {H_{\bf w}}{Z}\cong Z_2^2;$ if $p=5$ then
$\frac {H_{\bf w}}{Z}\cong A_4$; if $p=7$ then
$\frac {H_{\bf w}}{Z}\cong S_4$; and if $p=11$ then $\frac {H_{\bf w}}{Z}\cong A_5$
where $H_{\bf w}\cong \SL(2,5)$. Let $P$ be a Sylow $p$-subgroup of group
$H$. Then $H=H_{\bf w}P$ where $H_{\bf w}\cap P=1$ and $F=N\rtimes (PH_{\bf w})$. Now we may
identify ${\cal U}$ with vector space ${\mathbb V}=\mathbb{V}(2,p)$. Let $\a=(1,0)$. Then
$H_{\a}=\{ \left(
\begin{array}{cc}1 &0\\c &1\end{array}\right) \di c\in F_p\}\cong
Z _p$. Since for any $h\in H$, we have
$$((H_{\bf w})^h)_\a=(H_{\bf w})^h\cap H_\a=1\,\quad {\rm and}\,\quad
|{\a }^{(H_{\bf w})^h}|=p^2-1=|(H_{\bf w})^h|.$$
This implies that  acting on ${\cal U}$, $(H_{\bf w})^h$ fixes
0 and is transitive on $\mathbb{V}\setminus \{0\}$. Therefore, $p$ vertices
$\{ {\bf w}^h\di h\in H\}$ have the same neighborhood in $\si $, a
contradiction.

For $H=\GL(2,p)$ where $p=3, 5$, we have completely same argument as last paragraph and get  a contradiction again.
\qed

\subsection{Proof of Theorem~\ref{main2}.(2)}

The proof of Theorem~\ref{main2}.(2) consists of the following two lemmas.

\begin{lem}\label{4.1} Suppose that $\Aut(\si)$ acts primitively on  ${\cal W}$.   Then $p=3$.
\end{lem}
 \proof  Again set $F=\Aut(\si )$. By Theorem~\ref{main1}.(1), $F$ acts imprimitively on ${\cal U}$. Let $\mathfrak{U}$ be a   $p$-block system of $F$ on ${\cal U}$ with  the kernel $F_{(\mathfrak{U})}$. Then $F_{(\mathfrak{U})}\le (S_p)^p.$  Suppose that  $F$ is
primitive on ${\cal W}$. Then   $F_{(\mathfrak{U})}$ is transitive on
 ${\cal W}.$

Clearly,  $F$ is neither  a diagonal type or twisted wreath product type. So we only need to deal with three
cases separately: $F$ is almost simple type, product type or affine
type.

 \vskip 3mm
 (i) $F$ is almost simple type.

\vskip 3mm Let $T=\soc(F).$ Then $T$ is transitive on ${\cal W}$. Since
$T$ is the unique minimal normal subgroup of $F$,  it follows that
$T\le F_{(\mathfrak{U})}$, which implies  that  $T$ is transitive on each block
in $\mathfrak{U}$. Thus  $T$ has two faithful representations with respective  degree $p$ and $p^3$,
which is impossible.

\vskip 3mm
 (ii) $F$ is product type.

\vskip 3mm In this case,  $F=(M\times M\times M)\rtimes H$, where
$M$ is a primitive and irregular group of  degree $p$, where
$p\ge 5$. Since $p^3\di\di |F|$ and $F_{(\mathfrak{U})}$ is transitive on ${\cal W}$, we
have $p^3\di |F_{(\mathfrak{U})}|$ and so $p\nmid |F/F_{(\mathfrak{U})}|$, a contradiction.

 \vskip 3mm
 (iii) $F$ is affine type.

\vskip 3mm In this case, $F=N\rtimes H,$  where $N\cong Z_p^3$ and
$H$ is an irreducible subgroup  of $\GL(3,p).$ Clearly, $N\le F_{(\mathfrak{U})}$ and
thus $H$ must be transitive on $\mathfrak{U}$. Therefore, $H$ has a subgroup
$M$ of index $p$. Let $P$ be a Sylow $p$-subgroup of $H$. Suppose that there
exists an element $h$ of order $p$  in $H\cap F_{(\mathfrak{U})}$. Since $N\lg h\rg $ is a $p$-subgroup
in $F_{(\mathfrak{U})}$, it is abelian, and then  $[h, N]=1$, a contradiction,  noting $F$ is an affine group.
Therefore, $|P|=p$ and then $H=PM$.  Set $H_1=H\cap \SL(3,p)$. Noting $P\le \SL(3,p)$,  we get $H_1=PM_1$, where
$M_1=H_1\cap M.$  Set $Z=Z(\SL(3,p))$. Then $Z\cong Z_k$ for $k=(3,
p-1)$. Therefore, in $\PSL(3,p)$, $|\o{H_1}: \o{M_1}|=|\o{P}|=p$.
Since $\o{H_1}$ is an irreducible subgroup which has a subgroup
of index $p$, by checking Proposition~\ref{PSL3p}, the possible
candidates are
  $\PSL(2,5)$ or  $\PGL(2,5)$ for $p=5$;   $\PSL(2,7)$ for $p=7$;
   $\PSL(2,11)$ for $p=11$; and $Z_{13}\rtimes Z_3$, $A_4$ or $S_4$ for
   $p=3$.
 Moreover, if  $p=3,5,11$, then $H_1\cong \o{H_1}$; if $p=7$, the $H_1\cong  \o{H_1}$ or $H_1\cong  \o{H_1}\times Z_3.$
 In what follows, we shall show  $p\ne 5, 7, 11$ and then $p=3$, the lemma is proved.

\vskip 3mm

  For the contrary, suppose that  $p\in \{5,7,  11\}.$
 Since $H=PM$, we get that $H_{\mathfrak{u}}=M$  and $F_{\mathfrak{u}}=NM$, for some  $\mathfrak{u}\in \mathfrak{U}.$
 Then $\o{F_{\mathfrak{u}}}=F_{\mathfrak{u}}/F_{(\mathfrak{u})}
=\o{N}\o{M}\le S_p$. Since $\o{F_{\mathfrak{u}}}$ contains a normal regular
subgroup   $\o{N}$, it is an affine group, which implies that
$M/(M\cap F_{(\mathfrak{u})})\cong \o{M}\cong Z_l$ for $l\di (p-1).$ In particular,
$M_1/(M_1\cap F_{(\mathfrak{u})})\cong Z_{l'}$ for $l'\di (p-1)$. Note that our
group $M_1=A_4$ or $S_4$ for $p=5$; $S_4$ or $S_4\times Z_3$ for $p=7$;
 and   $A_5$ for $p=11$. In all the
cases,  three exists a subgroup $M_2\cong
A_4$  which is contained in $F_{(\mathfrak{u})}$, that is, $M_2$ fixes $\mathfrak{u}$
pointwise. For any ${\bf u} \in \mathfrak{u}$, we have that $N_{\bf u}\cong Z_p^2$ and $N_{\bf u}M_2$
fixes $\mathfrak{u}$ pointwise.

Now let's consider the subgroup $N_{\bf u}M_2$. Let $K_0$ be the kernel of
$M_2$ acting on $N_{\bf u}$ by conjugacy. Then $K_0$ fixes a 2-dimensional
subspace pointwise. It is easy to see that  the subgroup of  $\SL(3,p)$ fixing
a 2-dimensional subspace pointwise is isomorphic to $Z_p^2\rtimes
Z_{p-1}$.  Since $p\nmid |M_2|$, we know that  $K_0$ is cyclic. But $A_4$
contains only one  cyclic normal subgroup, that is  1, and thus
$K_0=1$ and then $M_2$ acts faithfully on $N_{\bf u}$, or equivalently,
$M_2\lesssim \GL(2,p).$ However, $\GL(2,p)$ does not contain any subgroup isomorphic to  $A_4$,
a contradiction.\qed

\begin{lem}\label{4.2}  $\Aut(\si)\cong S_3\wr S_3$, $\si \cong \si (3)$ or $\si(6)$;  and $\G\cong \G(9)$ or $\G(18)$ defined in Example ~\ref{defi1}.
\end{lem}
\demo Suppose $p=3$.  Continue the proof of (iii) in last paragraph.  Then  $|{\cal U}|=9$,    $F=N\rtimes H$, $H_1=H\cap \SL(3,p)$,
and  $H_1=Z_{13}\rtimes Z_3$ $A_4$ or $S_4$.
Clearly, $H_1\ne
Z_{13}\rtimes Z_3.$  Hence, we let  $A_4\le H_1\le H\le L$, where   $L$ is the same group in Lemma~\ref{L}, that is  the subgroup of $\GL(3,3)$ consisting of all those $3\times 3$ matrices with only one nonzero entry in each row and column and  actually, $L\cong Z_2^3\rtimes S_3$, of order $48$.

 \vskip 3mm

 (i) First, suppose that  $H=L$. Then $|F|=|N||H|=3^4 2^4.$
 Considering  the imprimitive action of $F$  on ${\cal U}$, we know that $F\lesssim S_3\wr S_3$.
Since $|S_3\wr S_3|=3^4 2^4$, we get $F\cong S_3\wr S_3$. In fact, $N\rtimes L$ is really isomorphic to $S_3\wr S_3.$

Since $F=N\times L$, an affine group, we may identify   ${\cal W}$ with the 3-dimensional space $\mathbb{V}=\mathbb{V}(3,3)$. Let  ${\bf w}$ be zero vector.  Then  $F_{\bf w}=L.$  Take a  vertex  ${\bf u}\in {\cal U}$. Then $N_{\bf u}=Z_p^2$, and $L_{\bf u}$ is a Sylow 2-subgroup of $L.$ Therefore, $F_{\bf u}=N_{\bf u}\rtimes L_{\bf u}.$   Consider $N_{\bf u}$ as a 2-dimensional subspace, $L_{\bf u}$ must preserve it.
As in Example~\ref{defi1}, set
$$\begin{array}{lll}{\mathbb{V}}_0&=&\{ (0, b, c)\di b, c\in \GF(3)\}, \,  {\mathbb{V}}_1=\{ (a, 0, c)\di a, c\in \GF(3)\},\\
 {\mathbb{V}}_2&=&\{ (a, b, 0)\di a, b\in \GF(3)\}. \end{array}$$
Without loss of  generality,
set $N_{\bf u}={\mathbb{V}}_0$.    Take an element $x={\footnotesize \left(\begin{array}{lll} 0 &1 &0\\  0 &0 &1 \\ 1 &0 &0\end{array}\right)}\in L.$
Then $\lg x\rg $ permutes ${\mathbb{V}}_0$, ${\mathbb{V}}_1$ and ${\mathbb{V}}_2$. Now, ${\cal U}$ may be identified with the set of nine lines:
$${\cal U}=\{ \a +{\mathbb{V}}_i\di \a \in {\mathbb{V}}, i\in Z_3\}.$$

It is easy to check that $F_{\bf w}(=L)$ has two orbits on ${\cal U}$ of length 3 and 6, respectively. Therefore, we just get two graphs, which are exactly $\si(3)$ and $\si(6)$, with $d({\bf w})=3$ and 6.

  Set $\si=\si (3)$ or $\si(6)$. From the argument of  last
section, we know that  there exist no graphs whose automorphism
group acts primitively on ${\cal U}$ and so  $\Aut(\si )$
acts primitively on ${\cal W}$ and imprimitively on ${\cal U}.$ Since   $S_3\wr
S_3$ is the maximal imprimitive group of degree 9,   $\Aut(\si )\le
S_3\wr S_3\cong F$, and then  $\Aut(\si )=S_3\wr S_3$.

Correspondingly, we get $\G\cong \G(9)$ or $\G(18)$.

\vskip 3mm (ii) Secondly, suppose that  $A_4\le H\lvertneqq L$. Then
$|\Aut(\si)|\lvertneqq 3^42^3$ and  $A_4\le F_{\bf w}$.
Consider the action of $A_4$ on ${\cal U}$. Clearly, each subgroup $Z_3$
of $A_4$ is transitive on ${\cal U}$ and the normal subgroup $Z_2^2$ of
$A_4$ fixes each block setwise. If $Z_2^2$ fixes pointwise in a
block, then $Z_2^2$ fixes ${\cal U}$ pointwise, which forces that $F$
acts unfaithfully on ${\cal U}$. Therefore, acting  in each block,
$Z_2^2$ fixes one vertex and exchange other two vertices, which
implies that $F_{\bf w}$ has two orbits on ${\cal U}$ with respective length 3
and 6. Hence, we get  two graphs as  same as in (i), that is $\si=\si(3)$ and  $\si(6)$, contradicting to $|\Aut(\si )|=3^42^4$.
 \qed

\subsection{Proof of Theorem~\ref{main2}.(3)}
Before proving Theorem~\ref{main2}.(3), we first prove two group theoretical results.

\begin{lem} \label{p-element} Let $\mathbb{V}=\mathbb{V}(3,p)$ and  $G=\GL(3,p).$ Take
$x=\footnotesize{\left(\begin{array}{ccc} 1 &2 &2 \\ 0&1 &2 \\
0&0 &1\end{array}\right)}\in G.$ Then
\begin{enumerate}
\item[{\rm (1)}] $x$ fixes setwise only one 1-dimensional subspace
$\lg \a\rg$ for
 $\a =(0, 0, 1)\in  \mathbb{V}$ and only one  2-dimensional subspace $\mathbb{S}'=\{ (0,
a_2, a_3)\di a_2, a_3\in \GF(p)\}$. \item[{\rm (2)}] For any
2-dimensional subspace $\mathbb{S}$ not including $\a$, we have
$\mathbb{S}^{x^i}\cap \mathbb{S}^{x^j}\cap \mathbb{S}^{x^k}=\{0\},$ where $i, j, k\in \GF(p)$
are distinct.
\end{enumerate}
\end{lem}
\demo (1) Checking directly.

(2) Let $\mathbb{S}$ be a 2-dimensional subspace  and $\a\not\in \mathbb{S}$.
 Suppose that $0\neq \b \in \mathbb{S}^{x^i}\cap \mathbb{S}^{x^j}\cap \mathbb{ S}^{x^k},$
  where $i, j, k$ are distinct. Then $\b^{x^{-i}}, \b^{x^{-j}}, \b^{x^{-k}}\in
 \mathbb{S}$. Since $x$ does not fix $\lg\b \rg $, the subspace
 $\lg \b^{x^{-i}}, \b^{x^{-j}},$
 $\b^{x^{-k}}\rg $ can  not be  1-dimensional and so it is $\mathbb{S}$.
  Set $\b= (a_1, a_2, a_3)$.
 Note that  for any $l\in F_p$,
 $$x^l={\footnotesize \left(\begin{array}{ccc} 1 &2l &2l^2\\ 0&1 &2l \\
0&0 &1\end{array}\right)}.$$ Let
$$D=\left| \begin{array}{c} \b^{x^{-i}}\\ \b^{x^{-j}} \\\b^{x^{-k}} \end{array} \right|
=\left| \begin{array}{ccc}
a_1 &-2i a_1+a_2 &2i^2a_1-2ia_2+a_3\\
a_1 &-2j a_1+a_2 &2j^2a_1-2ja_2+a_3\\
a_1 &-2k a_1+a_2 &2k^2a_1-2ka_2+a_3\\
 \end{array} \right| .$$
Since $\a \not\in \mathbb{S},$ we get $a_1\ne 0$. By computing we
get
$$D= 4a_1^3(i-j)(k-i)(k-j)\ne 0,$$
forcing dim($\mathbb{S}$)=3, a contradiction.
 \qed

\begin{lem}\label{shortlength}
Let $G$ be an imprimitive transitive group of degree $p^2$ on $\O$, where
$p\ge 3$ and $p^3\di |G|$ and  let $\mathcal{\mb } $ be  an imprimitive $p$-block system of $G$.
 Let $P$ be a Sylow p-subgroup of $G$. Then
\begin{enumerate}
\item[{\rm (1)}]  $\Exp(P)\le p^2$, $|Z(P)|=p$ and $P=(P\cap G_{(\mb)})\lg t\rg $, for some $t\in P$ such that
$t^p\in Z(P);$
\item[{\rm (2)}] Suppose that  provided either
$p=3$ or $p\ge 5$ and $ |P\cap G_{(\mb)}|\le p^{p-1}.$ Then $G_{(\mb)}$ is solvable, $P\cap G_{(\mb)}$ is a characteristic subgroup of  $G_{(\mb)}$ and so $P\cap G_{(\mb)}\lhd G$.
 \end{enumerate}
 \end{lem}
\demo (1) Check easily.

 \vskip 3mm
 (2) If $p=3$, then the conclusion  is clearly true.

 \vskip 3mm Suppose
$p\ge 5$ and $ |P\cap G_{(\mb)}|\le p^{p-1}.$
Set
$$K=G_{(\mb)}, \quad N=P\cap G_{(\mb)}, \quad \mb =\{ {\bf b}_0, {\bf b}_1, \cdots, {\bf b}_{p-1}\}.$$
  For any $g\in
K\setminus \{1\} $, let $\ell (g)$ be the number of blocks ${\bf b}_i$ in $\mb $ such
that the induced action $g^{{\bf b}_i}$ is nontrivial and set
$$\ell =min\{
\ell (g)\di g\in K\setminus \{ 1 \} \}.$$
Since $p^3\di |G|$ and $|N|\le p^{p-1}$, we get  $\ell
\ne p, 1$. Hence, $2\le \ell \le p-1$.

Take $g\in K$ such that $\ell(g)=\ell .$ Without loss of generality, say $g^{{\bf b}_i}$ is
nontrivial for $0\le i\le \ell -1$ and trivial for $\ell \le i\le
p-1$. Set $L=\cap_{\ell \le i\le {p-1}} K_{({\bf b}_i)}.$  Then $L\lhd K$.
Since $g\in L$, we get $L^{{\bf b}_i}$ is nontrivial for  $0\le i\le \ell -1$.
Since $L\lhd K$, it follows that $L$ is transitive  on each such  ${\bf b}_i.$
 By the definition of  $\ell $, we know that $L$ is faithful
 on ${\bf b}_i$ and so $N\cap L\cong Z_p$.

Take an element $x\in P\setminus N$  such that ${\bf b}_0^x={\bf b}_1$. Since
$\lg x\rg $ is transitive on ${\cal B}$, we have that $x$ cannot fix
setwise any proper subset of ${\cal B}.$ Therefore,
 $$1\le |\{ {\bf b}_0, {\bf b}_1, \cdots {\bf b}_{\ell-1}\}\cap
\{{\bf b}_0^x, {\bf b}_1^x, \cdots, {\bf b}_{\ell-1}^x\}|\le {\ell-1}.$$
 For any $h=h_0h_1\cdots h_{\ell-1}\in L,$ where $1\ne h_i\in S^{{\bf b}_i}$, we have that
   $h^x=h_0^xh_1^x\cdots h_{\ell-1}^x$. Noting that $1\ne h_0^x\in S^{{\bf b}_1},$
  we know that  both $L$ and
$L^x$ are nontrivial on ${\bf b}_1$. For any $h=h_0h_1\cdots h_{\ell-1}$ and
$h'=h_0'h_1'\cdots h_{\ell-1}'$ in $L,$  we have $\ell([h,(h')^x])\le
\ell-1.$  It follows the minimality of the value  $\ell$ that
$[h,(h')^x]=1$, equivalently,  $[L, L^x]=1$.  In
particular, $[L^{{\bf b}_1}, (L^x)^{{\bf b}_1}]=1$. Since $L^{{\bf b}_1}\le
C_{S^{{\bf b}_1}}((L^x)^{{\bf b}_1})$ and $(L^x)^{{\bf b}_1}$ is transitive on ${\bf b}_1$,
it follows from a well-known theorem in permutation group theory   that $L^{{\bf b}_1}$  is regular, that is, $L^{{\bf b}_1}\cong Z_p$.
Moreover, since $L$ is
faithful on each ${\bf b}_i$ for $0\le i\le \ell-1$ by the arguments in last paragraph, we get   $L\cong Z_p$. It has been proved that
$K^{{\bf b}_1}$ contains a regular normal subgroup $L^{{\bf b}_1}$, and so  $K^{{\bf b}_1}$ is
solvable. This in turn implies $K$ is solvable.

Let $N_1$ and $N_2$ be two Sylow $p$-subgroups of $K$.
Since $K^{{\bf b}_i}$ has the unique subgroup  $Z_p$ for
each block ${\bf b}_i$, we get
$[N_1,N_2]=1.$ Now $N_1N_2$ is a $p$-subgroup of $K$, which forces that
$N_1=N_2$. Therefore, $N$char $K$ and then $N\lhd G$, as desired.  \qed

\vskip 3mm

\f {\bf Proof of Theorem~\ref{main2}.(3):}
Suppose that $F=\Aut(\si)$ acts
imprimitively on  ${\cal W}$. By  Theorem~\ref{main2}.(1), $F$ also acts
imprimitively on  ${\cal U}$, with an imprimitive
complete $p-$block system    $\mathfrak{U}$.   Clearly, $F_{(\mathfrak{U})}\ne 1.$
Considering the imprimitive action of $F$ on ${\cal U}$,  we find that
$F\lesssim S_p\wr S_p=(S_p)^p\rtimes S_p,$ and $F_{(\mathfrak{U})}\lesssim (S_p)^p.$
Let $P\in \Syl_p(F).$ Then $P\lesssim (Z_p)^p\rtimes Z_p.$
Set $N=P\bigcap F_{(\mathfrak{U})}$. Then $N\lesssim (Z_p)^p$.

In what follows,
we divide our proof into two cases   depending on whether or
not $F_{(\mathfrak{U})}$ acts transitively on ${\cal W}$.

\vskip 3mm

(1) {\it $F_{(\mathfrak{U})}$ acts transitively on ${\cal W}$.}

\vskip 3mm

Suppose that $F_{(\mathfrak{U})}$ acts transitively on ${\cal W}$.
 Then $N$ is also  transitive on ${\cal W}$. Since   $N$ is abelian, $N$   acts
regularly on ${\cal W}$, that is $N\cong Z_p^3$ and then $|P|=p^4$.
  Take ${\bf w} \in {\cal W}$.
 Then  $P=N\rtimes P_{\bf w} \cong Z_p^3\rtimes Z_p.$
Considering  the action of $P$ on ${\cal U}$, for ${\bf u} \in {\cal U}$ we have
that  $P_{\bf u}=N_{\bf u}\cong Z_p^2.$

 By
Lemma~\ref{shortlength}, $F_{(\mathfrak{U})}$  is solvable and $N\lhd F$. Therefore,
$F$ is an affine group, that is  $F=N\rtimes F_{\bf w}$, where $N$ is
identified with
 the translation normal subgroup of $\AGL(3,p)$ and $F_{\bf w}$ with a
 reducible  subgroup of $\GL(3,p)$. That is the case (3.1) in Theorem~\ref{main2}.

\vskip 3mm

(2) {\it $F_{(\mathfrak{U})}$ acts intransitively on ${\cal W}$.}

 \vskip 3mm
 Suppose that $F_{(\mathfrak{U})}$ acts intransitively on ${\cal W}$. Since $
F/F_{(\mathfrak{U})}\le S_p$, we get $p\di\di |F/F_{(\mathfrak{U})}|$. Hence   $|N|\ge p^2$ and so $F_{(\mathfrak{U})}$
induces $p^2$-blocks on ${\cal W}$. Therefore, the first conclusion  of Theorem~\ref{main2}.(3.2) holds.

  Let ${\bf w}$ be any vertex in ${\cal W}$. Then we deal with two cases separately.

 \vskip 3mm
  (2.1) Suppose that ${\bf w}$ is exactly adjacent to two blocks in $\mathfrak{U}$. Then we are in case
  Theorem~\ref{main2}.(3.2.1).

  \vskip 3mm
  (2.2) Suppose that  ${\bf w}$ is adjacent to at least three blocks in $\mathfrak{U}.$  Then in what follows we shall prove the conclusions of Theorem~\ref{main2}.(3.2.2).

  Since  $F/F_{(\mathfrak{U})}\lesssim S_p$, it acts faithfully on both $\mathfrak{W}$ and  $\mathfrak{U}$.
  Thus,  we get $p\ne 3$ and so    $p\ge 5$. Since $|P|\ge p^3$ and $P$ acts faithfully on
 ${\cal U}$, it follows that $P$ is nonabelian.

Now we show $|N|=p^2$. For the contrary, suppose that  $|N|\ge p^3$. Let $N_1$  be a normal subgroup of $P$
such that $|N_1|=p^3$ and $N_1\le N.$ Let $x_0\in P\setminus N$ and
$P_1=N_1\lg x_0\rg $. Then by Lemma~\ref{shortlength}.(1),
$|Z(P)|=p$, $x_0^p\in Z(P)$ and then $|P_1|=p^4$. Clearly, for any ${\bf w}\in {\cal W}$, we have  that $|(P_1)_{\bf w}|=|(N_1)_{\bf w}|\ge p$;
 and for any ${\bf u}\in
{\cal U},$ we have that $(P_1)_{\bf u}=(N_1)_{\bf u}\cong Z_p^2$
and $N_1$ is transitive on every block in ${\cal U}$. As the same
reason as in (1), the conjugacy action of $x_0$ on
$N_1$ can be identified with the action of $x$ on $\mathbb{V}(3,p)$, where
$x$ is define as in Lemma~\ref{p-element}.
Suppose that ${\bf w}$ is adjacent to $p$ vertices  in a block  $\mathfrak{u}\in \mathfrak{U}$.
Since the edge-transitivity  of $\si$, we get that  ${\bf w}$  is adjacent
to $p$  vertices  in any block such that one of whose vertex  is adjacent to ${\bf w}$. Considering
the actions  $N_1$   on ${\cal W}$ and ${\cal U}$, we know that the vertices in
${\bf w}^{N_1}$ have  the same neighborhood, a contradiction. Therefore, if ${\bf w}$ is adjacent to a block $\mathfrak{u}\in \mathfrak{U}$,
 then $(N_1)_{\bf w}$ fixes pointwise $\mathfrak{u}$, equivalently
$(N_1)_{\bf w}\le (N_1)_{\bf u}$, otherwise,  ${\bf w}$ is adjacent to
$p$ vertices  in  $\mathfrak{u}$. By the hypothesis,   we assume that ${\bf w}$ is
adjacent to at least three blocks  $\mathfrak{u}^{x^i}, \mathfrak{u}^{x^j}$ and $\mathfrak{u}^{x^k}$,
where $i, j, k$ are distinct in  $Z_p$. Then  $(N_1)_{\bf w}\le
((N_1)_{\bf u})^{x^i}\cap ((N_1)_{\bf u})^{x^j}\cap ((N_1)_{\bf u})^{x^k}$.  Identifying
$(N_1)_{\bf w}, ((N_1)_{\bf u})^{x^i},  ((N_1)_{\bf u})^{x^j}$ and $((N_1)_{\bf u})^{x^k}$ with
the subspaces of $\mathbb{V}(3,p)$,  we get from Lemma~\ref{p-element} that
$(N_1)_{\bf w}=1$,  a contradiction.

Since $|N|=p^2$, we get   $|P|=p^3$.
 Since $p\ge 5$,   we get from Proposition~\ref{Dob} that
$P\lhd F,$ namely, $P$ is a normal subgroup of $F$ acting regularly on ${\cal W}$.
By Lemma~\ref{shortlength}, $F_{(\mathfrak{U})}$ is solvable.
 Moreover, since $Z_p\cong PF_{(\mathfrak{U})}/F_{(\mathfrak{U})}\lhd F/F_{(\mathfrak{U})}\le S_p$, it follows that $F/F_{(\mathfrak{U})}$ contains a normal regular subgroup  on $\mathfrak{U}$ and so it is affine group.   In other words,   $F/F_{(\mathfrak{U})}\cong Z_p\rtimes Z_r,$ where $r\di (p-1)$. \qed

\section{Examples of  graphs}

In this section,  by defining three bi-coset graphs  we show the existences of  the graphs in the three cases of Theorem~\ref{main2}.(3).

\vskip 3mm
\f {\it Graph $\si_1(p):$} For $p\ge 5$,  let $F=N\rtimes (\lg x\rg \rtimes H)\le \AGL(3,p)$ where $N$ is the translation
subgroup of $\AGL(3,p)$, where
$$ x={\footnotesize \left(\begin{array}{ccc} 1 &2 &2 \\ 0&1 &2 \\
0&0 &1\end{array}\right)},\quad H=\lg {\footnotesize \left(\begin{array}{ccc} s^2t^{-1} &0 &0 \\ 0&s &0 \\ 0&0 &t\end{array}\right)}
\di s, t\in \GF(p)^*\rg \le N_{\GL(3,p)}(\lg x\rg ).$$
   Let $N_0=\lg t_{(1, 0, 0)}, t_{(0, 1, 0)}\rg \le N.$ With the notation of Definition~\ref{defi2}, set
$$L=\lg x\rg H,\,\,  R=N_0H, \,\, D=RL.$$
Then $|F:L|=p^3$ and $|F:R|=p^2$.
Since $R$ and $L$ are nonmaximal subgroups of $F$,  we get that $F$ acts imprimitively on both $[F:L]$ and $[F:R]$.

Let $\si _1(p)={\bf B}(F;L,R,D)$, a double coset graph.  Since $|D|/|R|=|L|/|L\cap R|=p,$
the degree of any vertex in $[F:L]$ is $p$.
Moreover, one may easily see that there exist no two vertices in $[F:L]$ having the same neighborhood.
Now $F$ acts edge-transitively on $\si $.   Since the proof of $\Aut(\si)\cong F$ depends on
several lemmas and take a long argument, we do not try to  write it  in this paper but
shall put it   in our further paper.  This graph satisfies the condition  of  Theorem~\ref{main2}.(3.1).
Finally, let $\G_1(p)$ be the graph expanded from $\si _1(p)$.

\vskip 3mm
\f {\it Graph $\si_2(p):$} For any prime $p\ge 5$, let $$\s =(0, 1, \cdots , p-1),\quad  \t =(0)(1, -1)\cdots
(\frac{p-1}2, \frac{p+1}2)\in S_p.$$ Then $\lg \s, \t\rg \cong D_{2p}.$
Let $M\cong S_p$, $H\le M$ and $H\cong S_{p-1}$. Set
  $$\begin{array}{ll}&F=M\wr \lg \s, \t \rg=(\overbrace{M\times\cdots\times M}^{p\, {\rm times}})\rtimes \lg \s, \t\rg ,\\
&L=(\overbrace{M\times\cdots\times M}^{\frac{p-1}2-1\, {\rm times}}\times H\times H\times \overbrace{M\times\cdots\times M}^{\frac{p+1}2-1\, {\rm times}})\rtimes \lg \tau\rg,\\
&R=(H\times \overbrace{M\times\cdots \times M}^{p-1\, {\rm times}})\rtimes \lg \tau\rg , \,\, D=R\sigma ^{\frac{p-1}2}L.\end{array}$$
Then $|F:L|=p^3$ and  $|F:R|=p^2$.
 Clearly, $F$ acts imprimitively on both $[F:L]$ and $[F:R]$.

 Let $\si _2(p)={\bf B}(F;L,R,D)$.  Since $|D|/|R|=|L|/|L\cap R|=2,$
the degree of any vertex in $[G:L]$ is $2$.
    Moreover,     we shall  show  $\Aut(\si)\cong F$  in our further paper.
    Clearly,  $M^p$  induces a $p$-block system on $[F:R]$
and  the vertex $L$  is exactly adjacent to two blocks, corresponding to the double coset $R\sigma ^{\frac{p-1}2}L$.
This graph satisfies the condition  of  Theorem~\ref{main2}.(3.2.1).
Finally, let $\G_2(p)$ be the graph expanded from $\si _1(p)$.

\vskip 3mm
\f {\it Graph $\si_3(p):$} For any prime $p\ge 5$,  suppose that
$$P=\lg a, b\di a^{p^2}=b^p=1, [b, a]=c, a^p=c, a^b=a^{1+p}\rg .$$
Pick up  an element  $s$ of order $p-1$  in $Z_{p^2}^*.$ Let
$$F=P\rtimes \lg x\rg , \, {\rm where}\,   x^{p-1}=1, a^x=a^s$$
 Set
$$L=\lg x\rg ,\,\,  R=\lg b\rg \lg x\rg , \,\, D=RaL.$$
Then $|F:L|=p^3, |F:R|=p^2$.
 Clearly, $F$ acts imprimitively on both $[F:L]$ and $[F:R]$.

Let $\si _3(p)={\bf B}(F;L,R,D)$.  Since $|D|/|R|=|L|/|L\cap R|=p-1,$
the degree of any vertex in $[F:L]$ is $p-1$.
    Moreover,     we shall  show  $\Aut(\si)\cong F$  in our further paper.
    Clearly,  $\lg a^p\rg $  induces a $p$-block system on $[F:R]$
and  the vertex $L$  is adjacent to $p-1$ blocks. This graph satisfies the condition  of  Theorem~\ref{main2}.(3.2.2).
Finally, let $\G_3(p)$ be the graph expanded from $\si_3(p)$.

\vspace{3mm}
{\small
\f {\bf Acknowledgments:}  The  authors thank the referee for the helpful  comments and suggestions.
This work is partially   supported by
the National Natural Science Foundation of China and Natural Science Foundation of Beijing.}

\end{document}